\documentclass[11pt]{article}
\usepackage{graphicx} % Required for inserting images
\usepackage{amsmath,amssymb,enumerate,xcolor}

\title{Integrable geodesic flows  with simultaneously diagonalisable quadratic integrals}
\author{ Sergey I. Agafonov\footnote{Department of Mathematics,
S\~ao Paulo State University-UNESP, S\~ao Jos\'e do Rio Preto, Brazil,
 {\tt sergey.agafonov@gmail.com}}  \, and  Vladimir S.\ Matveev\footnote{
Institut f\"ur Mathematik, Friedrich Schiller Universit\"at Jena,
07737 Jena,  Germany  \ \ \quad {\tt  vladimir.matveev@uni-jena.de}}}

\newtheorem{theorem}{Theorem} 
\newtheorem{remark}{Remark} 
 
\newtheorem{corollary}{Corollary}[theorem]

\newcommand{\weg}[1]{}

\usepackage{biblatex} %Imports biblatex package
\begin{document}

\maketitle
\begin{abstract} 
We show that if $n$  functionally independent commutative quadratic in momenta integrals for the geodesic flow of a Riemannian  or pseudo-Riemannian metric on an $n$-dimensional manifold 
are simultaneously diagonalisable at the tangent space to every point, then they come from the St\"ackel construction, so the metric admits orthogonal separation of variables.  

{\bf MSC:  	37J35, 70H06}

{\bf Key words:} Quadratic in momenta integrals,  orthogonal separation of variables, finite-dimensional integrable systems,  Killing tensors, St\"ackel matrix
\end{abstract}

\section{Introduction} \label{sec:1}
We work locally, on a smooth $n$-dimensional pseudo-Riemannian manifold $(M,g)$ of any signature. By geodesic flow we understand the Hamiltonian system on $T^*M$ generated by the Hamiltonian 
$$H(x,p)= \tfrac{1}{2} g^{ij}p_ip_j. $$

We will study the situation when  the geodesic flow admits $n$, including the Hamiltonian,   integrals $$\overset{1}{I}(x,p)=2 H\, , \ \overset{2}{I}(x,p)\, , \ \dots \, , \  \overset{n}{I}(x,p)$$ 
such that the following conditions  are fulfilled: 

\begin{enumerate}
    \item The integrals are  quadratic in momenta, that is, $\overset{\alpha}{I}(x,p)=\overset{\alpha}{K}^{ij}p_ip_j$. In particular,  $g^{ij}=  \overset{1}{K}^{ij} $.  We  assume  without loss of generality that the (2,0)-tensor fields $\overset{\alpha}K^{ij}$ are symmetric in upper indexes. 
    \item At almost every point $x\in M$,  there exists a basis in $T_xM$ such that,  for every $\alpha=1,\dots, n$,   the matrix 
    $\left(\overset{\alpha}K^{ij}(x)\right)$ is diagonal.  
    \item  The differentials of the integrals are linearly independent at  least   one\footnote{Using ideas of \cite{Kruglikov_2016},  is is easy to show that  linear independence of the differentials of polynomial in momenta integrals  at one point implies their linear independence at almost every point, provided the manifold is connected}  point of $T^*M$.  
\end{enumerate}

In many publications on this topic, e.g. in \cite{BCR,eisenhart,Kiyohara},  it is assumed that for almost every point $x\in M$ the restrictions of  the tensor fields $\overset{\alpha}K^{ij}$, $\alpha=1,\dots, n$, 
to $T_xM$ are linearly  independent. Our main result, Theorem \ref{thm:1} below, shows that this assumptions follows from conditions (1,2,3): 
\begin{theorem} \label{thm:1}
 Under the assumptions  above,  for almost every point $x$   the  restrictions of the  tensor fields $\overset{\alpha}K^{ij}$, $\alpha=1,\dots, n$, 
to $T_xM$ are linearly  independent. In particular,   for a generic  linear combination  $I= \sum_{\alpha=2}^n \lambda_\alpha  \overset{\alpha}I$
 of the integrals,   the  corresponding 
 (1,1)-tensor field $K^i_j:= K^{si}g_{sj}$,  where   $K^{ij}= \sum_{\alpha=2}^n \lambda_\alpha  \overset{\alpha}{K}^{ij} $, has $n$ different eigenvalues. 

\end{theorem}

In dimension $n=3$, Theorem \ref{thm:1} and Corollary \ref{cor:1} below were proven in  \cite[Theorem 2]{Agafonov} by other methods.

\begin{corollary} \label{cor:1}
    Assume the integrals $\overset{\alpha}{I}$ satisfy the conditions (1,2,3) above and in addition are in involution with respect to the standard Poisson bracket. Then, near almost every point, the metric $g$   and the integrals   come from  the Stäckel construction.  
\end{corollary}
 
In view of Theorem \ref{thm:1},  Corollary  \ref{cor:1} follows from  \cite[Theorem 6]{KM}, 
\cite[Proposition 1.1.3]{Kiyohara}, \cite[Theorem 8.6]{BCR}  or, possibly,\footnote{By \cite[Note on page 185]{Klingenberg}  
the diploma thesis of A. Thimm 1976, which we were not able to find,  contains this result} from A. Thimm 1976. 
 In these references, it was shown that  $n$ quadratic functionally independent integrals in involution  such that the corresponding Killing tensors are  simultaneously diagonalisable at every tangent space and such that at least one of the Killing tensors with one index raised by the metric has $n$ different eigenvalues,  come from the St\"ackel construction which we recall below.    
 %Note that the paper \cite{eisenhart}  which is often cited in this context, is  not sufficient for deducing Corollary %\ref{cor:1} from Theorem \ref{thm:1}  as it is additionally assumed there that the metric and the integrals can locally be put in the diagonal form. By
% \cite[Theorem 8.6]{BCR},  this condition follows automatically from other conditions. 

As mentioned above,  the  difference between our conditions (1,2,3) and the assumptions  used  in \cite[Theorem 6]{KM}, \cite[Proposition 1.1.3]{Kiyohara} or  \cite[Theorem 8.6]{BCR} is as follows: in \cite[Theorem 6]{KM}, \cite[Proposition 1.1.3]{Kiyohara} or  \cite[Theorem 8.6]{BCR} it was assumed that one of the Killing tensors, with one index raised by the metric, has $n$ different eigenvalues. We do not have this condition as an assumption and  prove  that it follows from other assumptions.

Let us recall the St\"ackel\footnote{The construction appeared already in \cite[\S\S 13-14]{Liouville}, see also discussion in \cite[pp. 703--705]{Luetzen}} construction following  \cite{eisenhart,disser}. Take a  non-degenerate $n\times n$ matrix $S= (S_{ij})$  with $S_{ij}$ being a
function  of the $i$-th variable $x^i$ only . Next, consider the functions $\overset{\alpha}I$, $\alpha =1,\dots, n$,  given by the following system of linear equations
 \begin{equation} 
 \label{eq:St_intro}
S\mathbb{ I}  = \mathbb{P},  
\end{equation} 
where  $\mathbb{I}= \left(\overset{1}I, \overset{2}I,\dots ,\overset{n}I\right)^\top$ and $\mathbb{P}= \left(p_1^2, p_2^2,\dots ,p_n^2\right)^\top$.  It is known that the functions $\overset{\alpha}I$ are in involution. Taking one of them (say, the first one, provided the inverse matrix to $S$ has no zeros in the first raw) as twice 
the Hamiltonian of the metric, one obtains an integrable geodesic flow whose integrals satisfy the conditional (1,2,3). Corollary \ref{cor:1} says that locally, near almost every points,  there exist no other examples of geodesic flows admitting  $n$ independend quadratic in momenta 
integrals in involution,  such that the corresponding Killing tensors are simultaneously diagonalisable at almost every tangent space. 

It is known that metrics coming from the St\"ackel construction admit orthogonal separation of variables in the Hamilton-Jacobi equations, so the equation for their geodesics can be locally solved in quadratures \cite{BolsinovKonyaevMatveev-OrthogSep, book}. Namely,   J. Liouville \cite{Liouville} and, independently,  P. St\"ackel \cite{disser} has shown that  the metrics are precisely those admitting othogonal separation of variables. L.\,P.\, Eisenhart, in his widely cited and very influential  paper  \cite{eisenhart}, 
has shown that locally the  metrics coming from the St\"ackel construction are precisely those whose geodesic flows admit $n$ functionally independent integrals in involution satisfying the following conditions: the integrals are quadratic in momenta, 
the corresponding matrices  are simultaneously diagonalisable in a coordinate systems, and at every point the corresponding matrices are linearly independent. In \cite{BCR, KM, Kiyohara}  it was shown that the assumption that  the integrals  are  simultaneously diagonalisable in a coordinate system may be replaced  by a weaker assumption that the matrices of the integrals are diagonalisable in a frame. Our result further improves the result of Eisenhart and shows that the condition that the matrices of the integrals are linearly independent at each  point is not necessary as this assumption  follows from other conditions.

\subsubsection*{ Acknowledgements.} V.M.  thanks the DFG (projects 455806247 and  529233771) and  ARC Discovery Programme   (DP210100951) for the support. This research was partially supported by FAPESP (research fellowship 2018/20009-6), in particular the visit of V.M. to S\~{a}o Jos\'{e} do Rio Preto,  was  supported by this grant.  We thank  C. Chanu and G. Rastelli for pointing out  the references   \cite{BCR,KM} and for useful discussions.

\section{ Proof of Theorem \ref{thm:1}} \label{sec:2}
Under the assumptions (1,2,3) from Section \ref{sec:1}, 
near almost every point, there exists smooth  vector fields  
$v_1(x),...,v_n(x)\in T_xM$ such that  they are linearly independent at every tangent space and such that the metric $g$ and the
matrices $\overset{\alpha}K^{ij}$
 are diagonal in the basis $(v_1,...,v_n)$.  After re-arranging  and re-scaling the vectors $v_i$, there exists  $m \in \mathbb{N}$, $m\le n$, $k_1,...,k_m\in \mathbb{N}$ with 
 $k_1+\cdots + k_m=n$
 and  smooth local functions   $g_1 ,..,g_m , \overset{\alpha}\rho_1,\dots,\overset{\alpha}\rho_m,    $  \ $\alpha\in \{2,\dots, n\}$,  on $M$ such that 
 the Hamiltonian  and the integrals  $\overset{\alpha}I$  with  $ \alpha= 2,\dots, n$  
 are given by the formulas
 \begin{equation} \label{eq:2}
 \begin{array}{rcl}
 2H &=&   V_1 +  V_2+\cdots+   V_m\\
 \overset{\alpha}I &=& \overset{\alpha}\rho_1 V_1  + \overset{\alpha}\rho_2   V_2+\cdots+ \overset{\alpha}\rho_m  V_m \,  
 \end{array}
 \end{equation}
In the formulas above, $V_i$ are  the functions on the cotangent bundle given by 
$$\begin{array}{lcl} V_1&=&(v_1)^2\varepsilon_1+(v_2)^2\varepsilon_2+\cdots +(v_{k_1})^2\varepsilon_{k_1} ,  \\  V_2&=&(v_{k_1+1})^2\varepsilon_{k_1+1} +(v_{k_1+2})^2\varepsilon_{k_1+2} +\cdots +  (v_{k_1+k_2})^2\varepsilon_{ k_1+k_2}  ,  \\  & \vdots & \\ V_m&=& (v_{k_1+\cdots+k_{m-1}+1})^2\varepsilon_{k_1+\cdots+k_{m-1}+1}+   (v_{k_1+\cdots+k_{m-1}+2})^2\varepsilon_{k_1+\cdots+k_{m-1}+2}+\cdots +  (v_{n})^2\varepsilon_{n}, \end{array}$$
where $v_i$ is the linear function on the $T^*M$ generated by the vector field $v_i$ via the canonical 
identification\footnote{In naive terms, we consider the 
vector field $v=\sum_i v^i \partial_i $ as the linear function   $p\mapsto \sum_i v^ip_i$ on $T^*M$. This identification of vector fields on $M$ and linear in momenta functions on the cotangent bundle  is independent  of a  coordinate system} $TM\equiv T^{**}M$, and $\varepsilon_i\in \{-1, 1\}$.

The Poisson bracket of $H$ and $I= \overset{\alpha}I$ reads (we omit the index $\alpha$ since the equations hold for any $\overset{\alpha}I$): 
\begin{equation}\label{eq:1}
0=\{2 H, I\} =    \sum_{i,j=1}^m  \left( \{V_i,\rho_j\} V_j   + \rho_j \{V_i, V_j\}\right) .   
\end{equation}
The right hand side of  \eqref{eq:1} is a cubic polynomial in momenta so all its coefficients are zero. For every point $x\in M$, this gives us  
a system of linear equations on the directional derivatives $v_s(\rho_j)$  with $s\in \{1,\dots, n\}$ and  $j\in \{1,\dots, m\}$. The coefficients and free terms of this system  depend on $\rho_j(x)$, on the entries of the vector fields $v_s$ at $x$, and on the derivatives of the entries of the vector fields $v_s$ at $x$.    Let us show that 
all   directional derivatives  $v_s(\rho_i)$ can be reconstructed from this system.  
We will show this for the directional derivatives $v_1(\rho_2)$ and $v_1(\rho_1)$, since this will cover two principle cases $i=j$ and $i\ne j$;   the proof for all  other $v_i(\rho_j)$ is completely  analogous.  

In order to extract $v_1(\rho_2)$, note that  
the cubic in momenta component $  (v_{k_1+1})^2 v_1  $  shows up  only in the addends    
$\{V_1,\rho_2\} V_2$, $\rho_1 \{V_2, V_1\} $ and $\rho_1\{ V_1, V_2\}$. 
In these addends, the  coefficient  containing a derivative  of one of the  functions  $\rho_s$ is 
$v_1(\rho_2)$. Thus, equating the coefficient of    $   (v_{k_1+1})^2  v_1$  to zero   gives us $v_1(\rho_2)$ as a function of $\rho_1, \rho_2$ and the entries of $\{V_1, V_2\}$. 

Similarly, in order to extract $v_1(\rho_1)$,  we note that the cubic in momenta component $( v_1)^3  $  shows up  only in the addends    
$\{V_1,\rho_1\} V_1$  and $\rho_1\{ V_1, V_2\}$. 
Its coefficient containing the derivatives of $\rho$'s is 
$v_1(\rho_1)$. Thus, equating the coefficient of    $ (v_1)^3  $  to zero  gives us $v_1(\rho_1)$ 
as a function of $\rho_1$. 

Thus, all directional derivatives $v_s(\rho_j)$ can be obtained from the system \eqref{eq:1}. Let us now view the system \eqref{eq:1} as a linear  PDE-system on unknown functions  $\rho_i$.    The  coefficients of this system  come from the vector fields  
$v_s$ and are given by certain nonlinear expressions in the components of $v_s$ and their derivatives. 
Since the directional derivatives of all functions $\rho_i$ are expressed in the terms of  the functions 
$\rho_i$,    the system can be solved with respect to all derivatives of the functions $\rho_i$. Therefore, the initial  values of the functions $\rho_i$ at one point $x_0$ determine  the local solution of the system. This implies that the 
space of solutions is at most $m$-dimensional. Finally,  the linear vector space of the integrals $\overset{\alpha}I$ is at most $m$-dimensional. Since $n$ of them are functionally independent by our assumptions, $n=m$ and Theorem  \ref{thm:1} is proved.

\begin{remark}
    The  proof of Theorem \ref{thm:1} is  motivated by  \cite[proof of Lemma 1.2]{Benenti55-92}, \cite[\S 1.1]{Kiyohara} and  \cite[\S 2]{KKM}.  
\end{remark}

\printbibliography

@incollection {Benenti55-92,
    AUTHOR = {Benenti, S.},
     TITLE = {Inertia tensors and {S}t\"{a}ckel systems in the {E}uclidean
              spaces},
      NOTE = {Differential geometry (Turin, 1992)},
   JOURNAL = {Rend. Sem. Mat. Univ. Politec. Torino},
  FJOURNAL = {Universit\`a e Politecnico di Torino. Seminario Matematico.
              Rendiconti},
    VOLUME = {50},
      YEAR = {1992},
    NUMBER = {4},
     PAGES = {315--341 (1993)},
      ISSN = {0373-1243},
   MRCLASS = {70G10 (53C80 70H20)},
  MRNUMBER = {1261446},
MRREVIEWER = {Niky Kamran},
}

@article {BCR,
    AUTHOR = {Benenti, S. and Chanu, C. and Rastelli, G.},
     TITLE = {Remarks on the connection between the additive separation of
              the {H}amilton-{J}acobi equation and the multiplicative
              separation of the {S}chr\"{o}dinger equation. {I}. {T}he
              completeness and {R}obertson conditions},
   JOURNAL = {J. Math. Phys.},
  FJOURNAL = {Journal of Mathematical Physics},
    VOLUME = {43},
      YEAR = {2002},
    NUMBER = {11},
     PAGES = {5183--5222},
      ISSN = {0022-2488,1089-7658},
   MRCLASS = {35Q40 (81Q05)},
  MRNUMBER = {1935158},
MRREVIEWER = {A.\ S.\ Sumbatov},
       DOI = {10.1063/1.1506180},
       URL = {https://doi.org/10.1063/1.1506180},
}

@article {disser,
AUTHOR ={St\"ackel, P.},
TITLE = {Die Integration der Hamilton-Jacobischen Differentialgleichung mittelst Separation der Variablen},
JOURNAL = {Habilitationsschrift, Universit\"at Halle}, 
Year= {1891}
}

@article{BolsinovKonyaevMatveev-OrthogSep,
url = {https://arxiv.org/abs/2212.01605v2},
title = {Orthogonal separation of variables for spaces of constant curvature},
title = {},
  AUTHOR = {Bolsinov, A. V. and Konyaev, A. Yu. and  Matveev, V. S.},
journal = {Forum Mathematicum},
doi = {doi:10.1515/forum-2023-0300},
year = {2024},
lastchecked = {2024-03-18}
}

@article {Agafonov,
    AUTHOR = {Agafonov, S. I. },
     TITLE = {Integrable geodesic flow in 3D and webs of maximal rank}, 
JOURNAL = {arXiv:2403.01459},  YEAR = {2024},
}

@article{Liouville,
     author = {Liouville, J.},
     title = {M\'emoire sur l'int\'egration des \'equations diff\'erentielles du mouvement d'un nombre quelconque de points mat\'eriels},
     journal = {Journal de Math\'ematiques Pures et Appliqu\'ees},
     pages = {257--299},
     publisher = {Gauthier-Villars},
     volume = {1e s{\'e}rie, 14},
     year = {1849},
     language = {fr},
     url = {http://www.numdam.org/item/JMPA_1849_1_14__257_0/}
}

@book {Luetzen,
    AUTHOR = {L\"{u}tzen, J.},
     TITLE = {Joseph {L}iouville 1809--1882: master of pure and applied
              mathematics},
    SERIES = {Studies in the History of Mathematics and Physical Sciences},
    VOLUME = {15},
 PUBLISHER = {Springer-Verlag, New York},
      YEAR = {1990},
     PAGES = {xx+884},
      ISBN = {0-387-97180-7},
   MRCLASS = {01A70 (01A55)},
  MRNUMBER = {1066463},
MRREVIEWER = {F.\ Smithies},
       DOI = {10.1007/978-1-4612-0989-8},
       URL = {https://doi.org/10.1007/978-1-4612-0989-8},
}

@article {KM,
    AUTHOR = {Kalnins, E. G. and Miller,  W. jun.},
     TITLE = {Killing tensors and variable separation for
              {H}amilton-{J}acobi and {H}elmholtz equations},
   JOURNAL = {SIAM J. Math. Anal.},
  FJOURNAL = {SIAM Journal on Mathematical Analysis},
    VOLUME = {11},
      YEAR = {1980},
    NUMBER = {6},
     PAGES = {1011--1026},
      ISSN = {0036-1410},
   MRCLASS = {53B20 (35C99 70H20)},
  MRNUMBER = {595827},
MRREVIEWER = {A.\ P.\ Stone},
       DOI = {10.1137/0511089},
       URL = {https://doi.org/10.1137/0511089},
}

@Book{book,
 Author = {Kalnins, E. G. and Kress, J. M. and Miller, W. jun.},
 Title = {Separation of variables and superintegrability. {The} symmetry of solvable systems},
 FSeries = {IOP Expanding Physics},
 Series = {IOP Expand. Phys.},
 %ISBN = {978-0-7503-1315-5; 978-0-7503-1314-8},
 Year = {2018},
 Publisher = {Bristol: IOP Publishing},
 Language = {English},
 DOI = {10.1088/978-0-7503-1314-8},
 Keywords = {35-02,37Jxx,35C05},
 zbMATH = {7088725},
 Zbl = {1416.35002}
}

@article{Kiyohara,
  title={Two classes of Riemannian manifolds whose geodesic flows are integrable},
  author={ Kiyohara, K.},
  journal={Memoirs of the American Mathematical Society},
  year={1997},
  volume={130},
  pages={},
  url={https://api.semanticscholar.org/CorpusID:122681058}
}

@book {Klingenberg,
    AUTHOR = {Klingenberg, W.},
     TITLE = {Lectures on closed geodesics},
    SERIES = {Grundlehren der Mathematischen Wissenschaften},
    VOLUME = { 230},
 PUBLISHER = {Springer-Verlag, Berlin-New York},
      YEAR = {1978},
     PAGES = {x+227},
      ISBN = {3-540-08393-6},
   MRCLASS = {53C20 (58E10)},
  MRNUMBER = {478069},
MRREVIEWER = {Y.\ Mut\^{o}},
}

@article{KKM,
title = {When a (1,1)-tensor generates separation of variables of a certain metric},
journal = {Journal of Geometry and Physics},
volume = {195},
pages = {105031},
year = {2024},
issn = {0393-0440},
doi = {https://doi.org/10.1016/j.geomphys.2023.105031},
url = {https://www.sciencedirect.com/science/article/pii/S0393044023002838},
author = {Konyaev, A. Yu. and  Kress, J. M. and  Matveev, V. S.}
}

@article{Kruglikov_2016,
url = {https://dx.doi.org/10.1088/0951-7715/29/6/1755},
year = {2016},
publisher = {IOP Publishing},
volume = {29},
number = {6},
pages = {1755--1768},
author = { Kruglikov, B. S.  and   Matveev, V. S. },
title = {The geodesic flow of a generic metric does not admit nontrivial integrals polynomial in momenta},
journal = {Nonlinearity}
}

@article{eisenhart,
    AUTHOR = {Eisenhart, L. P.},
     TITLE = {Separable systems of {S}tackel},
   JOURNAL = {Ann. of Math. (2)},
  FJOURNAL = {Annals of Mathematics. Second Series},
    VOLUME = {35},
      YEAR = {1934},
    NUMBER = {2},
     PAGES = {284--305},
      ISSN = {0003-486X},
   MRCLASS = {DML},
  MRNUMBER = {1503163},
       DOI = {10.2307/1968433},
       URL = {https://doi.org/10.2307/1968433},
}

\end{document}